\documentclass[11pt]{article}
\usepackage{amsmath}
\usepackage{amsmath,amssymb,latexsym,color}
\usepackage[mathscr]{eucal}
\newtheorem{thm}{Theorem}[section]
\newtheorem{prop}[thm]{Proposition}
\newtheorem{lemma}[thm]{Lemma}

\newtheorem{example}{Example}[section]
\newtheorem{defin}{Definition}[section]

\newtheorem{remark}{Remark}[section]

\newcommand{\proof}{{\it Proof.\quad}}
\newcommand{\qed}{\hfill\Box\medskip}

\usepackage{CJK}
\begin{document}
\begin{CJK*}{GBK}{song}

\newcommand{\be}{\begin{equation}\label}
\newcommand{\ee}{\end{equation}}
\newcommand{\bea}{\begin{eqnarray}\label}
\newcommand{\eea}{\end{eqnarray}}

\title{\bf The smallest one-realization of a given set }

\author{
Ping Zhao$^{\rm a}$\quad Kefeng Diao$^{\rm a}$\quad    Kaishun
Wang$^{\rm b}$\thanks{Corresponding
author: wangks@bnu.edu.cn}\\
{\footnotesize a. \em School of Science, Linyi  University,
Linyi, Shandong, 276005, China }\\
{\footnotesize b. \em  Sch. Math. Sci. {\rm \&} Lab. Math. Com.
Sys., Beijing Normal University, Beijing 100875,  China} }

\date{}
 \maketitle

\begin{abstract}

For any set $S$ of positive integers, a mixed hypergraph ${\cal H}$
is a realization of $S$ if its feasible set is $S$, furthermore,
${\cal H}$ is a one-realization of $S$ if it is a realization of $S$
and each entry of its chromatic spectrum is either 0 or 1. Jiang et
al. \cite{Jiang} showed that the minimum number of vertices of
realization of $\{s,t\}$ with $2\leq s\leq t-2$ is $2t-s$.
Kr$\acute{\rm a}$l \cite{Kral} proved that there exists a
one-realization of $S$ with at most $|S|+2\max{S}-\min{S}$ vertices.
In this paper, we improve Kr$\acute{\rm a}$l's result, and determine
the size of the smallest one-realization of a given set. As a
result, we partially solve an open problem proposed by Jiang et al.
in 2002 and by Kr$\acute{\rm a}$l  in 2004.

\medskip
\noindent {\em Key words:} mixed hypergraph; feasible set; chromatic
spectrum; one-realization
\end{abstract}

\section{Introduction}

 A {\em mixed hypergraph } on a finite set $X$ is a triple ${\cal
H}=(X, {\cal C}, {\cal D})$, where ${\cal C}$ and ${\cal D}$ are
families of subsets of $X$, called the {\em ${\cal C}$-edges} and
{\em ${\cal D}$-edges}, respectively.  A \emph{bi-hypergraph} is a
mixed hypergraph with ${\cal C}={\cal D}$.
 A sub-hypergraph ${\cal H}'=(X', {\cal C}', {\cal D}')$ of a mixed
hypergraph ${\cal H}=(X, {\cal C}, {\cal D})$  is a {\em partial
sub-hypergraph}  if $X'=X$, and
 ${\cal H}'$ is called a {\em derived sub-hypergraph}
  of  ${\cal H}$ on $X'$, denoted by ${\cal H}[X']$, when ${\cal C}'=\{C\in {\cal C}| C\subseteq X'\}$
 and ${\cal D}'=\{D\in {\cal D}| D\subseteq X'\}$. Two mixed hypergraphs
${\cal H}_1=(X_1, {\cal C}_1, {\cal D}_1)$ and ${\cal H}_2=(X_2,
{\cal C}_2, {\cal D}_2)$ are \it isomorphic \rm if there exists a
bijection $\phi$ from $X_1$ to $X_2$ that maps each ${\cal C}$-edge
of ${\cal C}_1$ onto a ${\cal C}$-edge of ${\cal C}_2$ and maps each
${\cal D}$-edge of ${\cal D}_1$ onto
 a ${\cal D}$-edge of ${\cal D}_2$,
 and vice versa. The bijection $\phi$ is called an \emph{isomorphism} from ${\cal H}_1$ to ${\cal H}_2$.

 A {\em proper $k$-coloring} of
${\cal H}$ is a mapping from $X$ into a set of $k$ {\em colors} so
that each ${\cal C}$-edge has two vertices with a {\em Common} color
and each ${\cal D}$-edge has two vertices with {\em Distinct}
colors. A {\em strict $k$-coloring} is a proper $k$-coloring using
all of the $k$ colors, and a mixed hypergraph is {\em $k$-colorable}
if it has a strict $k$-coloring. The maximum (minimum) number of
colors in a strict coloring of ${\cal H}=(X, {\cal C}, {\cal D})$ is
the {\em upper chromatic number} $\overline \chi ({\cal H})$
 (resp. {\em lower chromatic number} $\chi ({\cal H})$) of ${\cal H}$.
The study of   the colorings of mixed hypergraphs has made a lot of
progress since its inception \cite{Voloshin1}.
 For more information, we would like refer readers to \cite{Tuza, Voloshin2}.

  The set of all the  values $k$ such that ${\cal H}$ has a strict $k$-coloring is called the {\em feasible set}
  of ${\cal H}$, denoted by ${\cal F}({\cal H})$. For each $k$, let $r_k$ denote
the number of {\em partitions} of the vertex set. Such partitions
are called {\em feasible partitions}.  The vector $R({\cal
H})=(r_1,r_2,\ldots,r_{\overline\chi})$ is called the {\em chromatic
spectrum} of ${\cal H}$. A mixed hypergraph has a {\em gap at $k$}
if its feasible set contains elements larger and smaller than $k$
but omits $k$. A {\em gap of size $g$} means $g$ consecutive gaps.
If some gaps occur, the feasible set and the chromatic spectrum of
${\cal H}$ are said to be {\em broken}, and if there are no gaps
then they are called {\em continuous} or {\em gap-free}. If $S$ is a
set of positive integers, we say that a mixed hypergraph ${\cal H}$
is a {\it realization} of $S$ if ${\cal F}({\cal H})=S$. A mixed
hypergraph ${\cal H}$ is a {\it one-realization} of $S$ if it is a
realization of $S$ and all the entries of the chromatic spectrum of
${\cal H}$ are either 0 or 1. This concept was firstly introduced by
Kr$\acute{\rm a}$l \cite{Kral}.

Bujt$\acute{\rm a}$s \cite{Bujtas} gave a necessary and sufficient
condition for a set $S$  to be the feasible set of an $r$-uniform
mixed hypergraph. Jiang et al. \cite{Jiang} proved that a set $S$ of
positive integers is a feasible set of a mixed hypergraph if and
only if $1\notin S$ or $S$ is an interval. They also discussed the
bound on the number of vertices of a mixed hypergraph with a gap, in
particular, the minimum number of vertices of realization  of
$\{s,t\}$ with $2\leq s\leq t-2$ is $2t-s$. Moreover, they mentioned
that the question of finding the minimum number of vertices in a
mixed hypergraph with feasible set $S$ of size at least 3 remains
open. In \cite{zdw}, we obtained an upper bound on the minimum
number of vertices of $3$-uniform bi-hypergraphs with a given
feasible set. Kr$\acute{\rm a}$l \cite{Kral} proved that there
exists a one-realization of $S$ with at most $|S|+2\max{S}-\min{S}$
vertices, and proposed the following  problem: What is the number of
 vertices of the smallest mixed hypergraph whose spectrum is
 equal to a given spectrum $(r_1,r_2,\ldots,r_m)$?

 In this paper, we determine the size of the smallest one-realization of a given set and obtain
 the following result:

 \begin{thm} For any integers $2\leq n_s<\cdots <n_2<n_1$, let $\delta(S)$ denote the minimum size
 of one-realizations of $S=\{n_1,n_2,\ldots,n_s\}$. Then
  $$\delta(S)=\left \{
\begin{array}{ll}2n_1-n_s, & \mbox{if~ $n_1>n_2+1,$}\\
2n_1-n_s-1,&\mbox{if ~$n_1=n_2+1.$}\\
\end{array}
\right. $$
\end{thm}

As a result, we partially solve the above open problem proposed  by
Jiang et al. and by Kr$\acute{\rm a}$l.

\section{Proof of Theorem 1.1}

In this section we always assume that $S=\{n_1,n_2,\ldots,n_s\}$ is
a set of integers with $2\leq n_s<\cdots <n_2<n_1$.  We first show
that the number $\delta(S)$ given in Theorem 1.1 is a lower bound on
the size of the smallest one-realization of $S$, then construct two
families of mixed hypergraphs which meet the bounds.

Jiang et al. \cite{Jiang} discussed the bound on the number of
vertices of a mixed hypergraph with a gap.

\begin{prop} {\rm(\cite[Theorem~3]{Jiang})} If ${\cal H}=(X, {\cal C}, {\cal D})$
is an $s$-colorable mixed hypergraph with a gap at $t-1$, then
$|X|\geq 2t-s$. For $2\leq s\leq t-2,$ this bound is sharp.
\end{prop}

\begin{lemma}
  $$\delta(S)\geq \left \{
\begin{array}{ll}2n_1-n_s,  & \mbox{if~ $n_1>n_2+1,$} \\
2n_1-n_s-1,  &\mbox{if~ $ n_1=n_2+1.$}\\
\end{array}
\right. $$
\end{lemma}

\proof Assume that ${\cal H}=(X, {\cal C}, {\cal D})$ is a
one-realization of $S$.

\smallskip
 {\em Case 1.} $n_1>n_2+1$.
Then, ${\cal H} $ has a gap at $n_1-1$. By Proposition 2.1, we have
 $\delta(S)\geq 2n_1-n_s$.

\smallskip{\em Case 2.} $n_1=n_2+1$.
Suppose  $|X|\leq 2n_1-(n_s+2)$.  For any strict $n_1$-coloring
$c_1=\{C_1,C_2,\ldots,C_{n_1}\}$ of ${\cal H}$, there exist at least
$n_s+2$ color classes of size one. Suppose $C_1=\{\alpha_1\},
C_2=\{\alpha_2\},\ldots, C_{n_s+2}=\{\alpha_{n_s+2}\}$. For any
strict $n_s$-coloring $c_s$ of ${\cal H}$, there are the following
two possible cases:

\smallskip {\em Case 2.1.} There exist three vertices in $\{\alpha_1,
\alpha_2, \ldots, \alpha_{n_s+2}\}$ which fall into a common color
class under $c_s$. Suppose $\alpha_1, \alpha_2, \alpha_3$ are in a
common color class under $c_s$. Then $\{\alpha_1, \alpha_2\},
\{\alpha_1, \alpha_3\}, \{\alpha_2, \alpha_3\}\notin {\cal D},$
 which follows that
 $\{C_1\cup
C_2,C_3, \ldots,C_{n_1}\}, \{C_1\cup C_3, C_2, C_4,
\ldots,C_{n_1}\}, \{C_1, C_2\cup C_3, C_4, \ldots,C_{n_1}\}$ are
strict $n_2$-colorings of ${\cal H}$. Therefore, ${\cal H}$ is not a
one-realization of $S$, a contradiction.

\smallskip

{\em Case 2.2.} There exist two pairs of vertices in $\{\alpha_1,
\alpha_2, \ldots, \alpha_{n_s+2}\}$
 each   of which falls into a common color class under $c_s$.
Suppose $\alpha_1, \alpha_2$ are in a common color class and
$\alpha_3, \alpha_4$ are in common color class under $c_s$. Then
$\{\alpha_1, \alpha_2\}, \{\alpha_3, \alpha_4\}\notin {\cal D},$ it
follows that $\{C_1\cup C_2,C_3, \ldots,C_{n_1}\}$ and $\{C_1,
C_2,C_3\cup C_4, C_5, \ldots,C_{n_1}\}$
 are strict $n_2$-colorings
of ${\cal H}$. Then  ${\cal H}$ is not a one-realization of $S$, a
contradiction. Hence, $\delta(S)\geq 2n_1-n_s-1$.$\qed$

In the rest of this section,  we shall construct two families of
mixed hypergraphs which meet the bound in Lemma 2.2.

For any positive integer $n$, let $[n]$ denote the set
$\{1,2,\ldots,n\}$.

\medskip

\noindent {\bf Construction I.} For any positive integer $s\geq 2,$
let
 \begin{eqnarray*}X_{n_1,\ldots,n_s}^0&=&\{(\underbrace{i,i,\ldots,i}_s)|~i=1,2,\ldots,n_s-1\},\\
X_{n_1,\ldots,n_s}^1&=&\bigcup_{t=2}^{s}\bigcup_{j=n_t}^{n_{t-1}-1}
\{(\underbrace{j,\ldots,j}_{t-1},n_t,n_{t+1},\ldots,n_s),
(\underbrace{j,\ldots,j}_{t-1},\underbrace{1,\ldots,1}_{s-t+1})\}.
 \end{eqnarray*}
Suppose
\begin{eqnarray*}&&X_{n_1,\ldots,n_s}^*=X_{n_1,\ldots,n_s}^0\cup X_{n_1,\ldots,n_s}^1\cup \{(n_1,n_2,\ldots,n_s)\},\\
&&{\cal
D}^*_{n_1,\ldots,n_s}=\{\{(x_1,x_2,\ldots,x_s),(y_1,y_2,\ldots,y_s)\}|x_i\neq y_i, i\in [s]\},\\
&&{\cal
C}^*_{n_1,\ldots,n_s}=\{\{(x_1,x_2,\ldots,x_s),(y_1,y_2,\ldots,y_s),(z_1,z_2,\ldots,z_s)\}|~
|\{x_j, y_j, z_j\}|=2, j\in [s]\}.
\end{eqnarray*}
Then ${\cal H}_{n_1,\ldots,n_s}^*= (X_{n_1,\ldots,n_s}^*,{\cal
C}_{n_1,\ldots,n_s}^*, {\cal D}_{n_1,\ldots,n_s}^*)$   is a mixed
hypergraph with  $2n_1-n_s$ vertices.

\medskip

Let
$$
\begin{array}{c}
X_{n_1,\ldots,n_s}=\{(x_1,x_2,\ldots,x_s)|x_i\in [n_i], i\in
[s]\},\\
 X^s_{ij}=\{(x_1,x_2,\ldots,x_{i-1},j, x_{i+1},\ldots,x_s)|x_k\in
[n_k], k\in [s]\setminus \{i\}\}, j\in [n_i].
 \end{array}
 $$
Then, for any $i\in [s]$,
$$c_i^{s*}=\{X_{i1}^*,X_{i2}^*,\ldots,X_{in_i}^*\}$$ is a strict
$n_i$-coloring of ${\cal H}_{n_1,\ldots,n_s}^*$, where
$X_{ij}^*=X_{n_1,\ldots,n_s}^*\cap X_{ij}^s, j\in [n_i]$.

\begin{lemma}  ${\cal H}_{n_1,n_2}^*$ is a one-realization of
$\{n_1, n_2\}$.
\end{lemma}

\proof Under any strict coloring $c=\{C_1,C_2,\ldots,C_m\}$ of
${\cal H}_{n_1,n_2}^*$, the vertices $(1,1), (2,2), \ldots,
(n_2,n_2)$ fall into distinct color classes. For each $i\in [n_2]$,
suppose $(i,i)\in C_i$. Then, for any $i\in [n_2-1]$ and $j\in
[n_1-n_2-1]$, we have $(n_2+j, n_2)\notin C_i$ and $(n_2+j,1)\notin
C_{n_2}$. Since $\{(1,1),(n_2,1),(n_2,n_2)\}$ is a ${\cal C}$-edge,
$(n_2,1)\in C_1$ or $C_{n_2}$.

\smallskip
{\em Case 1.} $(n_2,1)\in C_1$. The fact that
$\{(n_2,1),(n_2,n_2),(n_2+1,n_2)\}$ is a ${\cal C}$-edge follows
that  $(n_2+1,n_2)\in C_{n_2}$. From the ${\cal C}$-edge $
\{(n_2,1),(n_2+1,1),(n_2+1,n_2)\}$, we observe $(n_2+1,1)\in C_1$.
Similarly,
  $(n_2+j,1)\in C_1, (n_2+j,n_2)\in C_{n_2}$ for any $j\in [n_1-n_2-1]$ and $(n_1,n_2)\in C_{n_2}$.
  Therefore, $c=c_2^{2*}$.

\smallskip{\em Case 2.} $(n_2,1)\in C_{n_2}$.
The ${\cal D}$-edge $\{(n_2,1),(n_2+1,n_2)\}$ implies that
$(n_2+1,n_2)\notin C_{n_2}$. Suppose $(n_2+1, n_2)\in C_{n_2+1}$.
Owing to the ${\cal C}$-edge $\{(n_2,1),(n_2+1,1),(n_2+1,n_2)\}$, we
have  $(n_2+1,n_2)\in C_{n_2+1}$. Similarly,
$(n_2+j,n_2),(n_2+j,1)\in C_{n_2+j}$ for any $j\in [n_1-n_2-1]$ and
$(n_1,n_2)\in C_{n_1}$. Therefore, $c=c_1^{2*}$.

Hence, the desired result follows. $\qed$

\begin{thm}
${\cal H}_{n_1,\ldots,n_s}^*$ is a one-realization of $S$.
\end{thm}

\proof By Lemma 2.3,  the conclusion is true for $s=2$.

Let $X'=\{(x_2,x_2,x_3,x_4,\ldots,x_s)|x_j\in [n_j], j\in
[s]\setminus \{1\}\}$. Then ${\cal H}'={\cal
H}^*_{n_1,n_2,\ldots,n_s}[X']$ is isomorphic to ${\cal
H}_{n_2,n_3,n_4,\ldots,n_s}^*$. By induction, all the strict
colorings of ${\cal H}'$ are as follows:
$$c'_i=\{X'_{i1},X'_{i2},\ldots, X'_{in_i}\}, \quad i\in [s]\setminus
\{1\},$$ where $X'_{ij}=X'\cap X_{ij}^*, j\in [n_i]$.

For any strict coloring $c=\{C_1,C_2,\ldots,C_m\}$ of ${\cal
H}_{n_1,\ldots,n_s}^*$, the vertices $(1,1,\ldots,1),$
$(2,2,\ldots,2), \ldots, (n_s,n_s,\ldots,n_s)$  fall into distinct
color classes. Without loss of generality, suppose
$(i,i,\ldots,i)\in C_i$ for any $i\in [n_s]$. Then  there are the
following two possible cases:

 \smallskip
{\em Case 1.}  $c|_{X'}=c_2'$.
 The ${\cal C}$-edge $\{(1,1,\ldots,1),(n_2,1,\ldots,1),(n_2,n_2,n_3,\ldots,n_s)\}$
 implies that $(n_2,1,\ldots,1)\in C_1$ or $C_{n_2}$.

\smallskip
 {\em Case 1.1.} $(n_2,1,\ldots,1)\in C_1$. From the
${\cal C}$-edge $\{(n_2,1,\ldots,1),(n_2,n_2,n_3,\ldots,n_s)$,
$(n_2+1,n_2,n_3,\ldots,n_s)\}$
 and the ${\cal D}$-edge $\{(1,1,\ldots,1),(n_2+1,n_2,n_3,\ldots,n_s)\}$,
we observe $(n_2+1,n_2,n_3,\ldots,n_s)\in C_{n_2}$.
 By the ${\cal C}$-edge
$\{(n_2,n_2,n_3,\ldots,n_s),(n_2,1,\ldots,1),(n_2+1,1,\ldots,1)\}$
and the ${\cal D}$-edge
$\{(n_2,n_2,n_3,\ldots,n_s),(n_2+1,1,\ldots,1)\},$ we observe
$(n_2+1,1,\ldots,1)\in C_1$. Similarly, $(n_2+j,1,\ldots,1)\in C_1$,
 $(n_2+j,n_2,n_3,\ldots,n_s)\in C_{n_2}$ for any
 $j\in [n_1-n_2-1]$ and $(n_1,n_2,\ldots,n_s)\in C_{n_2}$. Therefore, $c=c_2^{s*}$.

\smallskip{\em Case 1.2.} $(n_2,1,\ldots,1)\in C_{n_2}$.
Note that $(n_2+j,1,\ldots,1)\notin C_k$ for any $j\in [n_1-n_2-1]$
and $k\in [n_2]\setminus \{1\}$. If $(n_2+1,1,\ldots,1)\in C_1$,
 from the ${\cal C}$-edge
$\{(n_2+1,1,\ldots,1),(n_2,n_2,n_3,\ldots,n_s),(n_2+1,n_2,\ldots,n_s)\},$
we observe $(n_2+1,n_2,\ldots,n_s)\in C_1$ or $C_{n_2}$, contrary to
the fact that both $\{(1,1,\ldots,1),(n_2+1,n_2,\ldots,n_s)\},
\{(n_2,1,\ldots,1)$, $(n_2+1,n_2,\ldots,n_s)\}$ are ${\cal
D}$-edges. Then, $(n_2+1,1,\ldots,1)\notin C_1$. Suppose
$(n_2+1,1,\ldots,1)\in C_{n_2+1}$.  The ${\cal C}$-edge
$\{(n_2+1,1,\ldots,1),(n_2+1,n_2,n_3,\ldots,n_s),(n_2,1,\ldots,1)\}$
   implies $(n_2+1,n_2,\ldots,n_s)\in C_{n_2+1}$. Similarly, $(n_2+j,1,\ldots,1),(n_2+j,n_2,\ldots,n_s)\in C_{n_2+j}$
   for any $j\in [n_1-n_2-1]$ and
   $(n_1,n_2,\ldots,n_s)\in C_{n_1}$. Therefore,  $c=c_1^{s*}$.

\smallskip
{\em Case 2.}  There exists a $k\in [s]\setminus \{1,2\}$ such that
$c|_{X'}=c_k'$. In this case, we observe
$(n_2,n_2,n_3,\ldots,n_k,\ldots,n_s)\in C_{n_k}$. For any $j\in
[n_1-n_2-1]$, the ${\cal D}$-edge
$\{(n_2+j,1,\ldots,1),(n_2,n_2,n_3,\ldots,n_k,\ldots,n_s)\}$ implies
that $(n_2+j,1,\ldots,1)\notin C_{n_k}$. From the ${\cal C}$-edge
$\{(1,1,\ldots,1),(n_2,n_2,n_3,\ldots,n_k,\ldots,n_s),(n_2,1,\ldots,1)\}$
and the ${\cal D}$-edge
$\{(n_k,\ldots,n_k,n_{k+1},\ldots,n_s),(n_2,1,\ldots,1)\},$ we
observe $(n_2,1,\ldots,1)\in C_1$. For any $j\in [n_1-n_2-1]$, the
${\cal C}$-edge
$\{(n_2+j,1,\ldots,1),(n_2,1,\ldots,1),(n_2,n_2,n_3,\ldots,n_k,\ldots,n_s)\}$
implies that $(n_2+j,1,\ldots,1)\in C_1$.

For any $j\in [n_1-n_2-1]$, since
$\{(1,1,\ldots,1),(n_2+j,n_2,\ldots,n_s)\}$ is a ${\cal D}$-edge,
$(n_2+j,n_2,\ldots,n_s)\notin C_1$. Moreover, the ${\cal C}$-edge
$\{(n_2,n_2,n_3,\ldots,n_s),(n_2+j,1,\ldots,1),(n_2+j,n_2,n_3,\ldots,n_s)\}$
implies that $(n_2+j,n_2,n_3,\ldots,n_s)\in C_{n_k}$ for any $j\in
[n_1-n_2-1]$. The fact that
$\{(n_1,n_2,n_3,\ldots,n_s),(n_2,1,\ldots,1),(n_2,n_2,n_3,\ldots,n_s)\}$
is a ${\cal C}$-edge follows that $(n_1,n_2,n_3,\ldots,n_s)\in
C_{n_k}$.  Hence, $c=c_k^{s*}$.

By the above discussion, the desired result follows. $\qed$

Next, we shall construct another family of mixed hypergraph. In this
case, we need to delete the vertex $(n_2,1,\ldots,1)$ from ${\cal
H}_{n_1,n_2,\ldots,n_s}^*$.

\medskip
\noindent {\bf Construction II.} Let
$X''=X_{n_1,n_2,\ldots,n_s}^*\setminus \{(n_2,1,\ldots,1)\}$ and
${\cal H}''={\cal H}_{n_1,n_2,\ldots,n_s}^*[X'']$. Then, for any
$i\in [s]$,
$$c_i''=\{X_{i1}'',X_{i2}'',\ldots,X_{in_i}''\}$$ is a strict
$n_i$-coloring of ${\cal H}''$, where $X_{ij}''=X''\cap X_{ij}^s,
j\in [n_i]$.

 \begin{thm} If $n_1=n_2+1$, the
${\cal H}''$ is a one-realization of $S$.
\end{thm}

\proof Referring  to the proof of Theorem 2.4,  all the strict
colorings of ${\cal H}_{n_2,n_2,n_3,\ldots,n_s}^*$ are
$$c'_i=\{X'_{i1},X'_{i2},\ldots, X'_{in_i}\}, \quad i\in [s]\setminus
\{1\},$$ where $X'=\{(x_2,x_2,x_3,x_4,\ldots,x_s)|x_j\in [n_j], j\in
[s]\setminus \{1\}\}$ and $X'_{ij}=X'\cap X_{ij}^*, j\in [n_i]$.

For any strict coloring $c=\{C_1,C_2,\ldots,C_m\}$ of ${\cal H}''$,
there are the following two possible cases:

\smallskip{\em Case 1.}  $c|_{X'}=c_2'$.
For any $(x_2,x_2,x_3,\ldots,x_s)\in X''$, suppose
$(x_2,x_2,x_3,\ldots,x_s)\in C_{x_2}$ under the coloring $c$.
 By the proof of Theorem 2.4, $(n_1,n_2,n_3,\ldots,n_s)\notin C_j$ for any $j\in [n_2-1]$.
Then, there are the following two possible subcases.

\smallskip
{\em Case 1.1.} $(n_1,n_2,n_3,\ldots,n_s)\in C_{n_2}.$
 It is immediate that $c=c_2''$.

  \smallskip
{\em Case 1.2} $(n_1,n_2,n_3,\ldots,n_s)\notin C_{n_2}.$ Suppose
$(n_1,n_2,n_3,\ldots,n_s)\in C_{n_1}$. Then, it is immediate that
 $c=c_1''$.

\smallskip{\em Case 2.}  There exists a $k\in [s]\setminus \{1,2\}$ such
that $c|_{X'}=c_k'$. It is immediate that
$(n_k,\ldots,n_k,n_{k+1},\ldots,n_s)\in C_{n_k}$ and
$(\underbrace{n_k,\ldots,n_k}_{k-1},1,\ldots,1)\in C_1$. From the
${\cal C}$-edge $\{(n_1,\ldots,n_k,
n_{k+1},\ldots,n_s),(n_k,\ldots,n_k,n_{k+1},\ldots,n_s),(n_k,\ldots,n_k,1,\ldots,1)\}$
and the ${\cal D}$-edge $\{(n_1,n_2,\ldots,n_s), (1,1,\ldots,1)\}$,
we observe $(n_1,n_2,\ldots,n_s)\in C_{n_k}$. Hence, $c=c_k''$.

Hence,  the desired result follows. $\qed$

Combining Lemma 2.2, Theorems 2.4 and  2.5, the proof of Theorem 1.1
is completed.

\section*{Acknowledgment}

The research
 is supported by NSF of
Shandong Province (No. ZR2009AM013), NCET-08-0052, NSF of China
(10871027) and the Fundamental Research Funds for the Central
Universities of China.

\end{CJK*}

\end{document}